\documentclass[11pt]{amsart} 

\usepackage{fullpage}

\usepackage{url}

\usepackage{amssymb}

\usepackage{cite}

\usepackage{graphicx}

\usepackage[usenames]{color}

\usepackage{accents}

\renewcommand{\a}{\alpha}
\renewcommand{\b}{\beta}
\newcommand{\g}{\gamma}
\renewcommand{\d}{\delta}
\newcommand{\D}{\Delta}
\newcommand{\e}{\varepsilon}
\newcommand{\f}{\varphi}

\newcommand{\Si}{\Sigma}

\renewcommand{\l}{\lambda}

\renewcommand{\O}{\Omega}

\newcommand{\cC}{{\mathcal C}}
\newcommand{\cM}{{\mathcal M}}

\newcommand{\cS}{{\mathcal S}}
\newcommand{\cB}{{\mathcal B}}
\newcommand{\cL}{{\mathcal L}}
\newcommand{\cE}{{\mathcal E}}
\newcommand{\cU}{{\mathcal U}}

\newcommand{\cP}{{\mathcal P}}

\newcommand{\cH}{{\mathcal H}}
\newcommand{\cD}{{\mathcal D}}

\newcommand{\cZ}{\mathcal Z}
\newcommand{\cI}{\mathcal I}

\newcommand{\bR}{\mathbb R}

\newcommand{\bE}{\mathbb E}

\newcommand{\bP}{\mathbb P}

\newcommand{\be}{\begin{equation}}
\newcommand{\ee}{\end{equation}}

\newcommand{\beaa}{\begin{eqnarray*}}
\newcommand{\bea}{\begin{eqnarray}}
\newcommand{\beal}[1]{\begin{eqnarray}\label{#1}}
\newcommand{\bean}{\begin{eqnarray}\nonumber}
\newcommand{\beadl}[1]{\begin{deqarr}\label{#1}}
\newcommand{\eeadl}[1]{\arrlabel{#1}\end{deqarr}}
\newcommand{\eeal}[1]{\label{#1}\end{eqnarray}}
\newcommand{\eead}[1]{\end{deqarr}}
\newcommand{\eea}{\end{eqnarray}}
\newcommand{\eeaa}{\end{eqnarray*}}

\newcommand{\Ric}{\operatorname{Ric}}

\renewcommand{\to}{\rightarrow}

\renewcommand{\phi}{\varphi}
\renewcommand{\epsilon}{\varepsilon}
\renewcommand{\hat}{\widehat}

\newcommand{\<}{\langle}
\renewcommand{\>}{\rangle}

\newcommand{\dm}{{\partial M}}

\newcommand{\w}{\widetilde}

\theoremstyle{plain}
\newtheorem{theorem}{Theorem}[section]

\newtheorem{remark}[theorem]{Remark}

\newtheorem{proposition}[theorem]{Proposition}

\newtheorem{conjecture}[theorem]{Conjecture}

\theoremstyle{definition}
\newtheorem{definition}[theorem]{Definition}

\def\endproof{\qed \medskip}
\def\blacksquare{\hbox to .60em {\vrule width .60em height .60em}}

\numberwithin{equation}{section}

\begin{document}

\title[ ]{Recent Progress and Problems on the Bartnik quasi-local mass}

\author{Michael T. Anderson}
\address{Dept. of Mathematics, 
Stony Brook University,
Stony Brook, NY 11790}
\email{anderson@math.sunysb.edu}

\thanks{Partially supported by NSF Grant DMS 1607479}

\begin{abstract}
This paper surveys recent progress on issues related to the Bartnik quasi-local mass $m_{B}$. In addition we formulate 
a number of new problems and conjectures regarding foundational properties of the mass $m_{B}$. This work is 
dedicated with pleasure to Robert Bartnik in honor of his $60^{\rm th}$ birthday.   

\end{abstract}

\maketitle

\setcounter{section}{0}
\setcounter{equation}{0}

\section{Introduction}
\setcounter{equation}{0}

  In this article, we survey recent progress on the Bartnik quasi-local mass $m_{B}$, discuss several problems 
with the various definitions and properties of $m_{B}$ and present new conjectures worthy of further study.

  We begin with a discussion of quasi-local mass in the time-symmetric or Riemannian setting of (space-like) 3-manifolds, 
which has been studied in much more detail.  
The constraint equations for $(g, K)$ then take the form $K = 0$ and $R_{g} \geq 0$ (dominant energy condition). In \S 5, 
we discuss the generalization to the space-time setting. (Of course any suitable definition of quasi-local mass should 
extend to the 4d space-time setting to be physically meaningful). 

  Let $\O$ be a compact, connected 3-manifold with connected boundary; typically $\O$ will be a closed $3$-ball $\bar B$ with 
$\partial \O = S^{2}$. Let $g_{\O}$ be a smooth metric of non-negative scalar curvature $R_{g_{\O}} \geq 0$ on $\O$. 
Let $\g = g_{\O}|_{\partial \O}$ be the induced metric on $\partial \O$ and let $A$ be 
the second fundamental form of $\partial \Omega$ in $(\O, g_{\O})$, with respect to the outward unit normal, 
so that $A  = \g$ for $S^{2}(1) \subset \bR^{3}$. 

  It is generally recognized that the basic properties desired of a quasi-local mass $m_{QL}(\O)$ of the domain $\O$ are 
the following: 

\smallskip 

\begin{itemize}

\item QL0: (Existence) $m_{QL}(\O)$ is well-defined on the class of domains $\O$ above. 

\smallskip

\item QL1: (Positivity) $m_{QL}(\O) \geq 0$ with equality if and only if $\O$ is flat, i.e.~there is an isometric embedding  
$F: (\O, g_{\O}) \to (\bR^{3}, g_{Eucl})$. 

\smallskip 

\item QL2: (Monotonicity) If $\O_{1} \subset \O_{2}$, then $m(\O_{1}) \leq m(\O_{2})$. 

\smallskip 

\item QL3: (Asymptotic behavior) If $\O_{i}$ is an exhaustion of a complete asymptotically flat manifold $(M, g)$ with 
$R_{g} \geq 0$, then $\lim m_{QL}(\O_{i}) = m_{ADM}$. 

\smallskip 

\item QL4: (Quasi-local) $m(\O)$ depends only on the first order geometry of $\partial \O$, i.e. $(\g, A)$. 

\end{itemize}

One could of course add further properties, such as the mass of a domain in the Schwarzschild metric surrounding the horizon 
be equal to the ADM mass $m$ of the Schwarzschild metric, small sphere limit behavior, etc.  

   Unfortunately, to the author's knowledge, there is no known definition of $m_{QL}(\O)$ which satisfies all of QL0-QL4, and 
perhaps there is no such quantity; we refer to [Sz] for a detailed discussion. 

\medskip 

   The Bartnik mass is a very natural and direct localization of the global ADM mass meant to capture the properties QL0-QL4. 
The original definition of Bartnik is the following, cf.~\cite{Ba1}, \cite{Ba2}. Let $\hat \cP^{0}$ denote the collection of smooth 
complete Riemannian 3-manifolds $(\hat M, \hat g)$ such that $\hat g$ is asymptotically flat (AF), of non-negative and integrable scalar 
curvature $R_{\hat g} \geq 0$, and $(\hat M, \hat g)$ has no horizons, i.e.~no (stable) minimal surfaces.  The Bartnik mass of a 
smooth bounded domain is then defined by 
\be \label{bmor}
m_{B}(\O) = \inf  \{m_{ADM}(\hat M, \hat g): (\O, g_{\O}) \subset (\hat M, \hat g) \in \hat \cP^{0}\},
\ee
where the infimum is taken over all $(\hat M, \hat g) \in \hat \cP^{0}$ for which $(\O, g_{\O})$ is isometrically embedded 
in $(\hat M, \hat g)$. 

  Many of the properties QL0-QL4 hold for $m_{B}$. Thus, Huisken-Illmanen \cite{HI} proved that $m_{B} \geq 0$ with 
equality if and only if $\O$ is locally flat, i.e.~there is an isometric immersion of $\O$ (or its universal cover) into $\bR^{3}$. 
(This cannot be improved to an embedding $\O \subset \bR^{3}$ by the work in \cite{AJ}). The property QL2 follows 
immediately from the definition and QL3 also follows from the work in \cite{HI}. 

  The main drawback of the definition \eqref{bm0} is that, in contrast to all other definitions of $m_{QL}$, it does not 
satisfy the quasi-local criterion QL4. This is because the no-horizon condition is global and depends in particular on 
the bounding domain $\O$, not just the geometry at $\partial \O$. Note that QL0 is also not fully satisfied, 
since $\O$ is not allowed to have any minimal surfaces. 

  Bartnik observed in \cite{Ba1}, \cite{Ba2} that an AF extension $(M, g)$ of $\partial \O$ with $\hat M = \Omega \cup M$ which 
realizes the infimum in \eqref{bmor} will in general only be Lipschitz along the ``seam'' $\partial \Omega = \dm$. By a 
simple and elegant argument using the $2^{\rm nd}$ variational formula for area, he showed that a minimizer should obey 
the boundary conditions 
\be \label{bcont}
\g_{\partial \O} = \g_{\dm}, \ \ \ \ H_{\partial \O} = H_{\dm},
\ee
where $H_{\dm}$ is the mean curvature of $\dm$ with respect to the unit normal pointing into $M$. The relation \eqref{bcont} 
implies that the scalar curvature is defined as a non-negative distribution across the seam. (Actually this holds for 
$H_{\partial \O} \geq H_{\dm}$ and this condition is sometimes used instead of \eqref{bcont}). Standard minimal surface 
arguments show that if $H_{\dm} \leq 0$ then any extension $(M, g)$ has a horizon, so that it is common practice to assume 
$$H = H_{\dm} > 0.$$

  The quasi-local data $(\g, H)$ on $\dm$ are now called Bartnik boundary data. Let $\cP_{(\g,H)}$ be the space of 
AF Riemannian manifolds $(M, g)$ with non-negative, integrable scalar curvature, inducing the data $(\g,H)$ at $\dm$. 
Let  $\cP_{(\g,H)}^{0} \subset \cP_{(\g,H)}$ be the subset of $(M, g)$ with no horizons, i.e.~immersed minimal $2$-spheres, 
surrounding $\dm$. A natural modification of the definition 
\eqref{bmor} which restores QL4 is then given by 
\be \label{bm}
m_{B}(\O) := m_{B}(\g, H) = \inf \{m_{ADM}(M, g): (M, g) \in \cP_{(\g,H)}^{0}\}. 
\ee

  The definition \eqref{bm} now satisfies QL4 and in addition QL0 holds for a larger class of domains $\O$. However, 
it remains an interesting open problem whether QL0 holds in general; given boundary data 
$(\g, H)$ (perhaps arising as boundary data of $(\O, g_{\O})$ with $R_{g_{\O}} \geq 0$), is it always true that 
\be \label{pne}
\cP_{(\g,H)} \neq \emptyset ?
\ee
This open problem is raised in \cite{Ba2}, cf.~also the recent results and survey given in \cite{AJ}. The same question holds 
with respect to $\cP_{(\g,H)}^{0}$, but here much less is known and the problem appears much harder. 
In addition, it is not clear or known whether the monotonicity property QL2 holds for the mass \eqref{bm}, 
again because of the no-horizon condition. 

   There have been a number of further modifications or variations of the definition of the Bartnik mass $m_{B}$; we 
refer to the recent paper by Jauregui \cite{J2} and further references therein for a careful and detailed discussion.  Further 
discussion of the no-horizon condition is given in \S 2 and \S 3. 

\medskip 

  It is to be expected that an extension $(M, g)$ realizing the infimum in \eqref{bm} or \eqref{bmor} satisfies strong 
conditions. In \cite{Ba1}, \cite{Ba2}, Bartnik presented a natural physical argument that extensions $(M, g)$ realizing the infimum in 
\eqref{bmor} (or \eqref{bm}) should be solutions of the static vacuum Einstein equations. Namely, any dynamical gravitational 
field carries energy and so mass (contributing to the Bondi mass at null infinity) and so an extension of minimal mass should 
have no gravitational dynamics, i.e.~be time-independent. For similar reasons, a minimal-mass extension should have no mass coming 
from matter sources, and so be vacuum. A time-independent vacuum solution which is time-symmetric ($K = 0$) is static vacuum. 

   The static vacuum Einstein equations are the equations for a pair $(g, u)$ on $M$ where $u$ is a potential function (forming 
the lapse function of the space-time $\cM = I\times M$) given by 
\be \label{stat}
uRic = D^{2}u, \ \ \D u = 0.
\ee
These equations are equivalent to the statement that the space-time metric $g_{\cM} = -u^{2}dt^{2} + g$ is Ricci-flat, i.e.~a vacuum 
solution of the Einstein equations. One usually adds the physical requirements that $u > 0$ and $u \to 1$ at infinity. A natural and 
physically well-motivated approach to proving this conjecture was developed in \cite{Ba3}, based on the Regge-Teitelboim Hamiltonian, 
cf.~also \cite{Ba4}. This is discussed further in \S 2. A full proof along these lines is given in \cite{AJ}. A completely different proof 
of part of this conjecture was given by Corvino in \cite{Co1} and more recently in \cite{Co2}. 

  There is a useful analogy of the Bartnik mass with the gravitational capacity of a body $\O \subset \bR^{3}$ in Newtonian gravity, 
or a charged body in electrostatics. Here one minimizes the Dirichlet energy,
\be \label{Eu}
E(v) = \int_{M}|dv|^{2},
\ee
over $M = \bR^{3}\setminus \O$ with boundary conditions $v = 1$ at $\partial \O$ and $v \to 0$ at infinity. (One could also 
set $v' = 1-v$ with $v' \to 1$ at infinity). Classical results show that the infimum of \eqref{Eu} is realized by a unique harmonic 
function $u$, $\D u = 0$ on $M$; $u$ represents the gravitational potential of the single layer $\partial \Omega$. The capacity 
of $\O$, equal to the total mass or charge up to a constant, is then given by 
\be \label{Nu}
E(u) = \inf E(v) = -\int_{\dm}N(u),
\ee
where $N$ is the unit normal into $M$ at $\dm$. 

   In \cite{Ba1}, \cite{Ba2}, Bartnik made the bold conjecture that a minimizer of \eqref{bmor} (or \eqref{bm}) also always 
exists and is unique; this is now called the minimization conjecture; 

\medskip 

{\bf Bartnik Minimization Conjecture.}  Any given boundary data $(\g, H)$ are realized by a unique AF minimizer $(M, g, u)$, 
$u > 0$, which is a solution of the static vacuum Einstein equations \eqref{stat} with boundary data $(\g,H)$. 

\medskip 
If true, this suggests the following conjecture \cite{Ba1}, \cite{Ba2}, which thus serves as a test 
of the minimization conjecture but is also of independent interest in geometric PDE theory. 

\medskip 

{\bf Bartnik Static Extension Conjecture.}  Given boundary data $(\g, H)$ on $\dm$, there exists a unique AF solution $(M, g, u)$, 
$u > 0$, of the static vacuum Einstein equations \eqref{stat} which induces the data $(\g, H)$ at $\dm$. 

\medskip 

  Unfortunately, the minimization conjecture has recently been proved to be false in general and both conjectures are very 
likely to be false in various regimes of boundary data $(\g, H)$. 
It thus becomes of basic interest to determine the realm of $(\g, H)$ for which these conjectures might be valid. 

  In the following sections of the paper, we present a discussion of various aspects of the Bartnik mass $m_{B}$ and 
and some of its further modifications.

\section{On the no-horizon condition}

  The no-horizon condition(s), i.e.~non-existence of minimal surfaces or minimal $2$-spheres discussed above, are a global 
requirement on the class of extensions over which one minimizes. Thus $m_{B}(\O)$ or $m_{B}(\g, H)$ depends on a large-scale, 
global property of the boundary data $(\g, H)$. The mass $m_{B}$ is undefined if $\cP_{(\g, H)}^{0} = \emptyset$, but it 
is not clear why. Moreover, it is very difficult to determine if a given AF extension with $R_{g} \geq 0$ has a minimal surface 
surrounding $\dm$ or no such surfaces. The same remarks hold for the global configuration $(\hat M, \hat g)$. This places 
a severe restriction on developing the foundations for a theory underlying the Bartnik mass or on understanding its behavior 
as a function of the boundary data $(\g,H)$. 

  The physical reasoning for a no-horizon condition is that mass behind a horizon is hidden (causally) from 
infinity and so is not captured by the ADM mass. One can always add enough matter sources (positive scalar curvature) 
to produce a horizon, so that the original body and such sources are hidden from infinity, where only an arbitrarily small mass 
may be detectable. Such configurations are the ``bag of gold" initial data sets coined by Wheeler \cite{W}, cf.~\cite{OM}. 
However, while the existence of a horizon allows for this behavior, the absence of horizons is not stable under time 
evolution. An initial time-symmetric slice $M$ (or $\hat M$) may have no horizons, and so be in $\cP_{(\g,H)}^{0}$ 
(or $\hat \cP^{0}$), but during the course of time evolution of the physical space-time horizons may develop. 

  This leads more naturally to the space-time Bartnik mass (cf.~\S 5). However, it is possible for instance that there exists a  
time-symmetric slice $\hat M'$ at a later time which has a horizon, so that such a slice is no longer in $\hat \cP^{0}$, even 
though the ADM mass remains conserved under such time evolution. This may occur even if the initial data for the 
Einstein evolution is vacuum, so $R_{g} = 0$. Of course restricting to the vacuum case, it is no longer so clear if the 
mass at infinity can be made arbitrarily small.  

  The main issue is to capture the gravitational mass. The mass due to matter sources is localized, in that one has 
an infinitesimal energy density $\rho: M \to \bR^{+}$. One can argue that these matter contributions should be subtracted in some 
sense from the total mass of an extension $(M, g)$ to obtain the gravitational mass of $(M, g)$. From this perspective, 
it again makes sense to restrict to vacuum configurations where $R_{g} = 0$. 

\medskip 

   Thus we propose the following modification of the Bartnik mass. Let $\cP$ be the space of smooth AF 
Riemannian $3$-manifolds $(M, g)$ with $R_{g} \geq 0$ and compact boundary $\dm$, and as above let 
$\cP_{(\g,H)} \subset \cP$ be the subset such that $(g|_{\dm}, H_{\dm}) = (\g, H)$; this is the space of extensions 
of the boundary data $(\g, H)$ with non-negative scalar curvature.  
Let 
$$\cC_{(\g,H)} = \{(M, g) \in \cP_{(\g,H)}: R_{g} = 0\},$$
be the set of solutions of the time-symmetric vacuum constraint equations ($R_{g} = 0$) with given Bartnik boundary data; 
(the constraint set with fixed boundary data).  

\begin{definition}
Let 
\be \label{bm0}
m_{B}^{0}(\g, H) = \inf \{m_{ADM}(M, g):  (M, g) \in \cC_{(\g,H)}\}.
\ee
\end{definition} 
A main point here is the absence of a horizon condition imposed on $(M, g)$. In particular $m_{B}^{0}(\g, H)$ 
satisfies the quasi-local criterion QL4. Some remarks related to this definition were given 
by Bartnik in the original work \cite{Ba1}, but the issue seems to have been neglected since then. 

  While the positive mass theorem (with corners) implies $m_{B}^{0}(\g, H) \geq 0$ if $(\g, H)$ are the boundary data 
of a compact domain $\O$ with non-negative scalar curvature,  the proof of Huisken-Illmanen \cite{HI} that $m_{B}^{0}(\O) > 0$ 
unless $\O$ is locally flat no longer applies. Thus we raise:

\begin{conjecture} 
One has 
$$m_{B}^{0}(\O) = 0,$$ 
only if $\O$ is locally flat. 
\end{conjecture} 
Of course if the infimum $m_{B}^{0}(\O)$ in \eqref{bm0} is realized by an extension of $\O$ in 
$\cC_{(g,H)}$, then the conjecture holds, by the positive mass theorem with corners.

\begin{remark} 
{\rm There is currently little evidence on which to base Conjecture 2.2. As suggested to the author by J.~Jauregui, consider 
a ``bag of gold" space $(M, g)$ with boundary data $(\g, H)$ and satisfying the following: $R_{g} \geq 0$, ${\rm supp}\, R_{g}$ 
compact in $M$ and $(M, g)$ equal to a Schwarzschild metric of arbitrarily small mass outside a compact set, as in the 
beginning of this section. Using results of Lohkamp \cite{L}, the metric $g$ may be deformed within a compact set to a new 
metric $\bar g$ with 
scalar curvature 
$$0 \leq R_{\bar g} < \e$$
and with $|\bar g - g|_{C^{0}} < \e$, for any prescribed $\e > 0$. The boundary data $(\g, H)$ and mass $m_{ADM}$ of $g$ remain 
the same for $\bar g$. Let $\w g = v^{4}\bar g$ where $v$ solves the equation 
$$-8 \w \D v + \w s v = 0,$$
with $v = 1$ on $\dm$ and $v \to 1$ at infinity. The metric $\w g$ is then scalar-flat $R_{\w g} = 0$ and satisfies 
$$0 \leq \w m < m, \ \  \w \g = \g, \ \ \w H \leq H,$$
(cf.~\cite{AJ}). Thus, one can produce a scalar-flat metric with arbitrarily small mass, and with smaller mean curvature at the 
boundary. This would lead to counterexamples to Conjecture 2.2 if (perhaps by some modified construction) $\w H$ could 
be made to equal $H$ (keeping $\w \g = \g$) or $\w H$ at least close to $H$ in a smooth topology. However, in the 
Lohkamp deformation above, the first and higher derivatives of $\bar g$ blow up as $\e \to 0$ and this may lead to blow-up of 
$\w H$ at $\dm$, so that $\w H$ is far away from $H$ in a smooth topology. It would be of 
interest to understand this situation better. 
  
}
\end{remark}

   The definition \eqref{bm0} has a natural variational interpretation closely related to work of Bartnik \cite{Ba4} in the 
case of complete manifolds without boundary, (cf.~discussion at the end of \cite{Ba3} and in \cite{Ba4}). This is 
based on the Regge-Teitelboim Hamiltonian and an approach to the positive mass theorem suggested by 
Brill-Deser-Fadeev \cite{BDF}.

 To describe this, let $\cS(M)$ denote the space of smooth pairs $(g, u)$ on $M$ where $g$ is an asymptotically flat metric and 
$u$ a function on $M$ with $u \to 1$ at infinity. Formally, the 4-metric $g_{\cM} = -u^{2}dt^{2} + g$ is a time-independent, 
i.e.~static, metric on the space-time $\cM = \bR\times M$. Of course $g_{\cM}$ is not vacuum in general. 
 
 Consider the Regge-Teitelboim Hamiltonian \cite{RT} in this setting:  
\be \label{cL}
\begin{array}{cc}
\cH_{RT}: \cS(M) \to \bR, \\
\cH_{RT}(g, u) = 16\pi m_{ADM}(g) - \int_{M}uR_{g} dv_{g}.
\end{array}
\ee
The first term differs from the Einstein--Hilbert action on the 4-manifold $\cM$ by a divergence term. 
In contrast to the original ADM Hamiltonian, the Regge--Teitelboim Hamiltonian \eqref{cL} (suitably regularized) 
gives a well-defined variational problem and is a smooth functional on the full manifold $\cS(M)$, cf.~\cite{RT}, \cite{Ba3} 
for details. 

  The first variation of the RT Hamiltonian is well-known, cf.~\cite{RT}, \cite{Ba2}, \cite{AJ}; 
\be \label{grad}
\nabla \cH = -(S^{*}u + {\tfrac{1}{2}}uR_{g} g, R_{g}, uA - N(u)\g, 2u)
\ee
in the sense that, if $(h, u')$ is any variation of $(g, u)$ inducing the variation $(h^{T}, H'_{h})$ of boundary data $(\g, H)$, 
then
\be \label{var}
-d\cH_{(g, u)}(h, u', h^{T}, H'_{h}) = \int_{M}[\<S^{*}u + {\tfrac{1}{2}}uR_{g}g, h\> + R_{g}u']dv_{g} + \int_{\dm}
[\<uA - N(u)\g, h^{T}\> + 2uH'_{h}]dv_{\g}.
\ee
Here $N$ is the unit normal at $\dm$ pointing into $M$, $A$ is the $2^{\rm nd}$ 
fundamental form of $\dm$ in $M$ and  
\be \label{S*}
S^{*}u = D^{2}u - (\D u) g - u\Ric,
\ee
is the formal $L^2$-adjoint of the linearization $R'=DR_g$ of the scalar curvature. The static vacuum Einstein equations 
\eqref{stat} are equivalent to the system $(g, u) \in \cS(M)$ such that 
\be \label{stat1}
S^{*}u = 0, \ \ R = 0.
\ee
(The condition $R = 0$ is equivalent to $\D u = 0$ when $S^{*}u = 0$). 

  Let $\cS_{(\g,H)}(M)\subset \cS(M)$ be the space of static metrics with fixed Bartnik boundary data $(\g,H)$. 
Working in a given function space such as a weighted H\"older space $C_{\d}^{m,\a}$, it is straightforward 
(using the implicit function theorem) to show that $\cS_{(\g,H)}(M)$ is a smooth, closed Banach submanifold of 
$\cS(M)$, for all choices of $(\g, H)$.

  Thus, \eqref{var} shows that critical points of the Hamiltonian $\cH_{RT}$ on $\cS_{(\g,H)}(M)$ are given by static vacuum 
Einstein metrics realizing the given boundary data $(\g, H)$. Note that although $u \to 1$ at infinity, it is not assumed at this point 
that $u > 0$ on $M$. 
\medskip 

  Now for the restriction of $\cH_{RT}$ to $\cC_{(\g,H)}$, one clearly has 
\be \label{Hm}
\cH_{RT} |_{\cC_{(\g,H)}} = 16\pi m_{ADM}(g): \cC_{(\g,H)} \to \bR.
\ee
It is proved in \cite{AJ} that, if non-empty, $\cC_{(\g,H)}$ is a smooth, infinite dimensional Banach manifold in a natural 
$C_{\d}^{m,\a}$ H\"older space topology, and that $m_{ADM}$ in \eqref{Hm} is a smooth functional on $\cC_{(\g,H)}$. 
Moreover, it is shown that critical points of $\cH_{RT}$ or $m_{ADM}$ on the constraint manifold $\cC_{(\g,H)}$ are exactly 
the static vacuum Einstein metrics $(g, u)$ satisfying the elliptic boundary conditions $(\g, H)$. In addition, minimizers of the 
Bartnik mass \eqref{bmor}, or \eqref{bm}, are static vacuum solutions with $u > 0$ and $u \to 1$ at infinity, cf.~\cite{AJ}. 

  It is worth noting that Miao \cite{Mi} has proved, based on the black hole uniqueness theorem, that static vacuum solutions 
$(M, g, u)$ as in \eqref{stat1} with $u > 0$ on $M$ have no minimal surfaces surrounding $\dm$. On the other hand, it is 
unknown whether static vacuum solutions $(M, g, u)$ have no compact minimal surfaces at all. 

  In contrast to the discussion above, there are no critical points of the mass 
$$m_{ADM}: \cS_{(\g,H)}(M) \to \bR,$$
so that some constraint, such as $R_{g} \geq 0$, is necessary to obtain a constrained critical point or minimizer of $m_{ADM}$. 
In fact if $R_{g}$ is not identically zero, i.e.~$g \notin \cC_{(\g,H)}$, then the mass $m_{ADM}$ can be 
decreased, both infinitesimally and locally, in the space $\cS_{(\g, H^{\leq})}^{+}(M)$, consisting of metrics of 
non-negative scalar curvature in $\cS(M)$ with mean curvature $\leq H$ and induced metric $\g$ at $\dm$, cf.~\cite{AJ}. 

\medskip 

  Next we relate the domains of definition of the masses in \eqref{bm} and \eqref{bm0}. 

\begin{conjecture}
If $(M, g)$ is an extension of the boundary data $(\g,H)$ in $\cP$, then there exists an extension in 
$\cC_{(\g,H)}$, i.e.
$$\cP_{(\g,H)} \neq \emptyset \Rightarrow \cC_{(\g,H)} \neq \emptyset.$$
\end{conjecture}
This is of course a very special case of (and much simpler than) the Bartnik minimization conjecture. 
It is not clear that either $\cP_{(\g,H)}$ or $\cC_{(\g,H)}$ are non-empty, for any given $(\g, H)$ arising as boundary 
data of a compact body $\O$ with $R_{g_{\O}} \geq 0$; a good testing ground would be the class of immersed spheres 
in $\bR^{3}$ discussed in \S 4. 

  We also note the following monotonicity property. Suppose $(\O_{1}, g_{\O_{1}}) \subset (\O_{2}, g_{\O_{2}})$ and 
$R_{g_{\O_{2}}} = 0$, so there are no matter sources in $\O_{2}$. Then 
\be \label{mono}
m_{B}^{0}(\O_{1}) \leq m_{B}^{0}(\O_{2}).
\ee
This follows immediately from the definition \eqref{bm0}, if one allows the space $\cC_{(\g,H)}$ to have Lipschitz seams 
across smooth surfaces. The condition $R_{g} = 0$ is then well-defined distributionally.  If equality holds in \eqref{mono} 
and $m_{B}^{0}(\O_{1})$ is realized by a static vacuum extension $(M_{1}, g_{1}, u_{1})$, then $(\O_{2}\setminus \O_{1}, g_{2}) 
\subset (M_{1}, g_{1})$ isometrically. 

\medskip 

  Next we discuss issues related to the Bartnik static extension conjecture. Via the discussion above, these are closely 
related to the study of $m_{B}$ itself. 

  Let $\bE_{st}$ be the set of all AF solutions $(M, g, u)$ of the static vacuum Einstein equations \eqref{stat}, $u > 0$ on $M$ and 
let $\cE_{st}$ be the associated moduli space, obtained from $\bE_{st}$ by quotienting out by the group of diffeomorphisms 
equal to the identity on $\dm$. Again, given a suitable H\"older space $C_{\d}^{m,\a}$ topology, it is proved in \cite{A1}, 
\cite{AK} that $\cE_{st}$ is a smooth infinite dimensional Banach manifold. Let $\cB$ be the space of Bartnik boundary data 
$(\g, H)$; again this may be given a natural smooth Banach manifold topology. Moreover, the Bartnik boundary map 
\be \label{bmap}
\begin{array}{cc}
\Pi_{B}: \cE_{st} \to \cB, \\
\Pi_{B}([g, u]) = (\g, H),
\end{array}
\ee
is a smooth Fredholm map, of Fredholm index zero. The Bartnik static extension conjecture from \S 1 
is thus the statement that the smooth map $\Pi_{B}$ is a bijection, i.e.~both surjective (existence) and injective (uniqueness). The 
formulation \eqref{bmap} gives an effective means to study the static extension conjecture in detail. 

  Note that the ADM mass
\be \label{mass}
m_{ADM}: \cE_{st} \to \bR,
\ee
is a smooth functional on $\cE_{st}$. It can be ``effectively" computed as a usual Komar integral, i.e.~
$$m_{ADM} = \int_{\dm}N(u);$$
compare with \eqref{Nu}. 
In regions where $\Pi_{B}$ is injective with smooth inverse, the mass $m_{B}(\g, H)$ or $m_{B}^{0}(\g, H)$ thus becomes 
a smooth functional of the boundary data $(\g, H)$, assuming it is realized as in the Minimization Conjecture. In regions 
where $\Pi_{B}$ fails to be injective, this is unlikely to be the case. 

   These remarks lead naturally to consideration of a further choice $m_{B}^{st}$ for the definition of Bartnik mass. Namely 
$m_{B}^{st}(\g, H)$ is defined only for 
$$(\g, H) \in {\rm Im}\, \Pi_{B}$$
and 
\be \label{bmstat}
m_{B}^{st}(\g,H) = \inf \{m_{ADM}(M, g, u): (M, g, u) \in \cE_{st} \ {\rm and} \ \Pi_{B}(M, g, u) = (\g,H)\}.
\ee
In many situations, one would expect that $(\Pi_{B})^{-1}(\g,H)$ is compact, and generically a finite number of points; this 
will be the case in regions $\cU \subset \cE_{st}$ where $\Pi_{B}$ is proper. In such cases, the infimum in \eqref{bmstat} 
is then realized as the smallest value of a continuous function on a compact set. Thus $m_{B}^{st}$ is, in principle, 
effectively computable. 

\medskip 

   Note that the mass functional \eqref{mass} does not map into $\bR^{+}$ in general, i.e.~it is not assumed in \eqref{mass} or 
\eqref{bmap} that the boundary data $(\g, H)$ bound a body $\O$ with $R_{g} \geq 0$. This is closely related to the 
{\it fill-in problem}, dual to the extension problem \eqref{pne} discussed above: given 
data $(\g, H)$ on $\partial \O$, when does there exist a metric $g_{\O}$ on $\O$ with $R_{g_{\O}} \geq 0$ 
and with Bartnik boundary data $(\g, H)$?  

  It is proved by Jauregui in \cite{J1} that for any $(\g, H)$, with Gauss curvature $K_{\g} > 0$ and $H > 0$, one has a basic 
trichotomy: there exists $\l_{0} > 0$ so that $(\g, H)$ has a fill-in with $R > 0$ for all $\l < \l_{0}$ and no fill-in with 
$R \geq 0$ for $\l > \l_{0}$. It is conjectured in \cite{J1} that the data $(\g, \l_{0}H)$ have a fill-in which is a solution 
to the static vacuum equations \eqref{stat}; thus static vacuum solutions provide a division between these two behaviors. 
We refer to \cite{A2} for an analysis of the behavior of static vacuum solutions on bounded domains $\O$. 

   Regarding the dual extension problem, it is proved in \cite{A3} that given any fixed $(\g, H)$ as 
above (with $\g$ now arbitrary), there is a $\l_{0}$ such that $(\g, \l_{0}H)$ are boundary data of a complete 
exterior AF static vacuum solution $(M, g, u)$, $M = \bR^{3}\setminus B$. Moreover, for $\l > \l_{0}$, the data 
$(\g, \l H)$ have an AF extension with non-negative scalar curvature, i.e.~$\cP_{(\g, \l H)} \neq \emptyset$, cf.~\cite{AJ}. 
On the other hand, the trichotomy above breaks down for $\l < \l_{0}$: Mantoulidis-Schoen \cite{MS} have constructed 
large families of data $(\g, \l H)$ with $\l > 0$ arbitrarily small such that $\cP_{(\g, \l H)} \neq \emptyset$. 

\medskip 

   Unfortunately, it is proved in \cite{AK} that $\Pi_{B}$ in \eqref{bmap} is not a proper map and so is not a homeomorphism, for a 
rather large domain within $\cE_{st}$ or $\cB$.  
At this time, it seems very unlikely that $\Pi_{B}$ is either globally surjective or globally injective, and so the domain in 
$\cE_{st}$ on which $\Pi_{B}$ is a smooth diffeomorphism onto its image becomes of basic importance; these issues are 
discussed further in \S 4.

\section{The outer-minimizing modification}

   An interesting modification of the Bartnik mass $m_{B}$ has been suggested by Bray \cite{BC}, based on considerations from the 
Riemannian Penrose inequality \cite{HI}, \cite{Br}. Given $(M, g) \in \cP$, 
$\dm$ is said to be outer-minimizing in $(M, g)$ if 
$$area (\dm) < area (\Si)$$
for any surface $\Si \subset M$ surrounding $\dm$, (so that the cycle $\dm - \Si$ is null-homologous). Let $\cP^{out}\subset \cP$ 
be the subset for which $\dm$ is outer-minimizing, so that $\cP^{out}$ is an open domain in $\cP$. Given Bartnik boundary 
data $(\g, H)$, let $\cP_{(\g,H)}^{out} = \cP_{(\g, H)} \cap \cP^{out}$ be the space of outer-minimizing extensions of $(\g, H)$. 
Define then 
\be \label{bmout}
m_{B}^{out}(\g,H) = \inf \{m_{ADM}(M, g): (M, g) \in \cP_{(\g,H)}^{out})\}.
\ee   

  The mass $m_{B}^{out}$ has a number of advantageous properties which follow from the fundamental work of Huisken-Illmanen 
\cite{HI}. Namely the Hawking mass $m_{H}$ is a lower bound for $m_{B}^{out}$: 
$$m_{H} \leq m_{B}^{out}.$$  
For domains $\O$ in the Schwarzschild metric $(M, g_{Sch})$ containing the Schwarzschild horizon $r = 2m$, one has 
\be \label{Sch}
m_{B}^{out}(\O) = m_{Sch}.
\ee
The relation \eqref{Sch} is unknown for any other version \eqref{bmor}, \eqref{bm} or \eqref{bm0} of the Bartnik mass. 
Directly from the definition, one sees that a restricted version of the monotonicity property QL2 also holds; 
if $\Omega_{1}\subset \O_{2}$ and $\O_{1}$ is outer-minimizing in $\O_{2}$ then 
$m_{B}^{out}(\O_{1}) \leq m_{B}^{out}(\O_{2})$. Further, QL3 holds. 

\medskip 

  Perhaps the main issue in extending the understanding of $m_{B}^{out}$ is the domain of existence property QL0. Namely, 
just as with the no-horizon condition, the definition $m_{B}^{out}$ again depends on a large-scale, global property of the 
boundary data $(\g, H)$. It is very difficult to effectively determine whether Bartnik boundary data $(\g, H)$ admit an 
outer-minimizing extension in $\cP$, i.e.~for which $(\g, H)$ is 
$$\cP_{(\g,H)}^{out} \neq \emptyset.$$
The mass $m_{B}^{out}$ is not defined when $\cP_{(\g,H)}^{out} = \emptyset$ but there is no solid physical reason for this. 

  Consider the simplest case of flat boundary data, 
i.e.~boundary data $(\g, H)$ of a flat bounded domain in $\bR^{3}$ for which $m_{B} = 0$. At least when $\partial \O = 
S^{2}$, flat boundary data $(\g, H)$ uniquely determine the bounding domain $\O$ (up to rigid motion). However there are 
very few general results on characterizing the data $(\g, H)$ for which the boundary is outer-minimizing in $\bR^{3}$. The wonderful 
solution of the Weyl embedding problem by Nirenberg \cite{N} and Pogorelov \cite{P} gives an affirmative answer when the Gauss curvature 
$K_{\g} > 0$; convex boundaries are outer-minimizing in $\bR^{3}$. However, there do not appear to be any sound physical 
reasons why the mass of a non-convex body in $\bR^{3}$ should not be well-defined; one would expect it is defined and 
the mass is zero. Thus it is of basic interest to understand more general conditions on $(\g, H)$ for which $\partial \O$ is 
outer-minimizing in $\bR^{3}$. Moreover, there there are no known analogs of the Nirenberg-Pogorelov result for 
general solutions $(M, g, u)$ of the static vacuum Einstein equations. 

\medskip 

  A further issue with $m_{B}^{out}$ is the failure of the Static Extension and Minimization conjectures in this context, 
i.e.~the non-existence of static vacuum extensions realizing $m_{B}^{out}$. This behavior occurs for a very natural and 
simple class of boundary data $(\g, H)$. The following result is essentially given in \cite{AJ}: 

\begin{proposition} Let $\g$ be a smooth metric on $S^{2}$ distinct from the round metric $\g_{2m}$ of radius $2m$. 
Then there is a constant $H_{m} > 0$, depending only on $m$ and a lower bound on the $C^{0}$-distance between $\g$ 
and $\{\g_{2m'}, m' \in \bR^{+}\}$, such that, for any smooth function $H$ satisfying
$$0 < H < H_{m},$$
pointwise, the Bartnik boundary data $(\g, H)$ are not the boundary data of a solution $(M, g, u)$ of the static vacuum Einstein 
\eqref{stat1} with $\dm$ outer-minimizing in $(M, g)$. 
\end{proposition}

{\bf Proof:} The proof is by contradiction. Thus, let $(\g_{i}, H_{i})$, $H_{i} > 0$, be a sequence of smooth Bartnik boundary data 
converging smoothly to $(\g, 0)$. Suppose for each $i$, $(\g_{i}, H_{i})$ are boundary data of static vacuum solutions 
$(M, g_{i}, u_{i})$, $u_{i} > 0$, with $\dm$ outer-minimizing in $(M, g_{i})$. By the basic compactness theorem in \cite{AK}, a 
subsequence of $(M, g_{i}, u_{i})$ converges, smoothly up to $\dm$, to a limit solution $(M, g, u)$, $u > 0$, of the static 
vacuum Einstein equations with Bartnik boundary data 
$$\lim_{i\to \infty} (\g_{i}, H_{i}) = (\g, 0).$$
This compactness theorem strongly uses (a weak form of) the outer-minimizing property. 

  By the black hole uniqueness theorem and its extension by Miao in \cite{Mi}, the only static vacuum solution $(M, g, u)$ with 
$u > 0$, $u \to 1$ at infinity and horizon $H = 0$ boundary is the Schwarzschild metric with boundary data $(\g_{2m}, 0)$ for some 
$m$. It follows that if $\g \neq \g_{2m'}$ for some $m'$, i.e.~$\g$ is not a round metric, then such a sequence cannot exist. 
Thus, there is a neighborhood of the boundary data $(\g, 0)$ where there are no outer-minimizing solutions of the 
static vacuum Einstein equations. 

{\endproof}

  In sum, the black hole uniqueness theorem causes the non-existence of static vacuum solutions $(M, g, u)$ with 
sufficiently small mean curvature and non-round boundary metric. It follows in particular that the mass $m_{B}^{out}$ 
cannot be realized by boundary data $(\g, H)$ with $H$ sufficiently small and $\g$ not a round metric. 

  We believe that Proposition 3.1 holds for all static vacuum solutions, without the outer-minimizing assumption. 
It follows from the work in \cite{AK} that if $(M, g, u)$ is a static vacuum solution with $\g$ not a round metric and $H$ sufficiently 
small, then the curvature of $(M, g)$ becomes either very large near $\dm$ or the injectivity radius normal to $\dm$ in $(M, g)$ becomes 
very small (or both). There does seem to any particular reason which would cause such behaviour with such controlled boundary data; 
however this remains an open problem.

\section{Embeddings and Immersions}

   It is proved in \cite{AJ} that there exists a large family of smooth immersions $F: S^{2} \to \bR^{3}$, which extend to 
smooth embeddings of the interior 3-ball $F: B^{3} \to \bR^{3}$, for which the corresponding Bartnik data $(\g, H)$ 
induced from $F$ have 
$$m_{B}(\g, H) = 0,$$
but $(\g, H)$ are not realizable as boundary data of mass-minimizing extensions. In particular, for such 
data $\cP_{(\g,H)}^{0} \neq \emptyset$ so that $m_{B}$ is defined but is not achieved by a mass-minimizer. 
This gives the first counterexamples to the Bartnik minimization conjecture; the discussion at the end of \S 3 suggests 
there are more. 

   This behavior illustrates that the manifold-with-boundary structure breaks down for all minimizing sequences 
for the Bartnik mass, even with smooth control of the boundary data $(\g, H)$. 

  We conjecture that a similar principle holds for general immersed 2-spheres. Let ${\rm Imm_{A}}(S^{2}, \bR^{3})$ be the space 
of smooth Alexandrov immersed $2$-spheres in $\bR^{3}$, i.e.~immersions of spheres which extend to immersions of the 3-ball 
$B^{3}$. The inverse image $F^{-1}(\bR^{3})$ is thus a locally flat 3-ball $(\O, g_{\O})$ and $F$ is the developing map giving an 
isometric immersion of the abstract body $\O \simeq B^{3}$ into $\bR^{3}$. Let $\cI = {\rm Imm_{A}}(S^{2}, \bR^{3}) \setminus 
{\rm Emb}(S^{2}, \bR^{3})$ be the complement of the space of embeddings $S^{2} \to \bR^{3}$. The constructions in \cite{AJ} 
take place at the boundary $\partial \cI$. For $F \in \cI$, let $(\g, H)$ be the induced metric and mean curvature from the 
Euclidean metric. 

\begin{conjecture} 
For $F \in \cI$ as above, either 
$$\cP_{(\g, H)} = \emptyset,$$
(so that $m_{B}(\g, H)$ is not defined), or $m_{B}(\g, H)$ is not achieved by a static vacuum solution with boundary data $(\g, H)$. 
\end{conjecture} 

  This behavior is not restricted to flat $\bR^{3}$; one may conjecture the same behavior for immersions of spheres into any 
given static vacuum solution $(M, g, u)$. 

\medskip 

   For $F \in \cI$ as above, the abstract body $(\O, g_{\O})$ becomes ''twisted" and has self-intersections when viewed as 
a body in the static vacuum solution $(\bR^{3}, g_{Eucl}, 1)$ (or more generally $(M, g, u)$). However, note that $(\O, g_{\O})$ 
may be smoothly approximated by embedded space-like 3-balls $(\O_{i}, g_{i}, K_{i})$ in Minkowski space time $\bR^{1,3}$, 
i.e. the 4-dimensional (vacuum) space-time development of the initial Cauchy hypersurface $(\bR^{3}, g_{Eucl}, K = 0)$. 
Thus there are always space-like AF extensions $(M, g, K)$ of $(\O, g_{\O}, 0)$ with approximately (and possibly even the same) 
boundary data in the 4-dimensional space-time; cf.~the discussion in \S 5. This passage from a fixed 3-dimensional 
space-like slice to the the full 4-dimensional space-time bears resemblance to the passage from the Brown-York mass \cite{BY} 
to the more recent Wang-Yau mass \cite{WY}. 
 
  In the case of $\bR^{3}$, the reasoning above suggests that the mass of any body $(\O, g_{\O})$ arising from 
$\cI$ should still be defined and be zero. As in the discussion in \S 3, there do not seem to be any physical reasons why the mass of an 
embedded body in $\bR^{3}$ should be defined but become undefined in passing continuously to an immersed body. 
Moreover, it is very difficult to detect from the boundary data $(\g, H)$ of a flat body $\O$ whether $\O$ gives an 
embedded or immersed sphere in $\bR^{3}$. 

\medskip 

 From this perspective, one may consider the following problem. Define an AF static vacuum solution $(M, g, u)$, $u > 0$, 
$M$ open, without boundary, to be {\it maximal} if it cannot be smoothly extended to any larger domain $(M', g', u')$, $u' > 0$, 
with $(M, g, u)$ isometrically embedded in $(M', g', u')$ as a proper subset; thus $M$ is maximal with respect to inclusion. Of 
course $(\bR^{3}, g_{Eucl}, 1)$ is a maximal solution; it is the only maximal solution which is complete as a metric space. 

\medskip 

  {\bf Problem.} Given Bartnik boundary data $(\g, H)$, is there a maximal static vacuum solution $(M, g, u)$ and an immersion 
$F: S^{2} \to M$ which realizes $(\g, H)$, so that the ``boundary data" on $F(S^{2})$ are $(\g, H)$. 

  Note that an affirmative resolution of the problem would then allow a definition of Bartnik mass for any $(\g, H)$; namely one 
may define $\w m_{B}(\g, H)$ to be the infimum of $m_{ADM}$ over all static vacuum solutions $(M, g, u)$ for which there is 
an immersion $F: S^{2} \to M$ inducing the data $(\g, H)$. 

\medskip 
 
  From one perspective, this problem may be viewed as an analog of the classical isometric immersion problem for 
spheres in $\bR^{3}$: given a metric $\g$ on $S^{2}$, is there an isometric immersion $F: S^{2} \to \bR^{3}$ realizing $\g$, 
so that $F^{*}(g_{Eucl}) = \g$? This is a famous and notoriously hard problem; the main result to date is the solution to the 
Weyl embedding problem by Nirenberg and Pogorelov. This generalization to the static vacuum case $(M, g, u)$ adds on the 
prescription of the mean curvature $H$, corresponding to the extra scalar field $u$. Note however that a positive solution of 
the static vacuum Problem above does not imply a resolution of the isometric immersion problem. One significant advantage of 
the static vacuum problem is that it is given by an elliptic system of equations, in strong contrast to the isometric 
immersion problem. Thus it might well be more amenable to further progress. 

   Observe however that the static vacuum Problem above is incompatible with the discussion following Proposition 3.1. 
Thus, a more reasonable version of the Problem is the following: 

\begin{conjecture} 
  Given Bartnik boundary data $(\g, H)$, with (say) $area_{\g}(S^{2}) = 1$, there is a constant $H_{1} > 0$ 
such that if 
$$H \geq H_{1} > 0,$$
then there is a maximal static vacuum solution $(M, g, u)$ and an immersion $F: S^{2} \to M$ which realizes $(\g, H)$, 
i.e.~the ``boundary data" induced on $F(S^{2})$ are $(\g, H)$. 
\end{conjecture} 

  Many of the problems raised here for the Bartnik mass exist also for other QL mass definitions. For instance, the 
Brown-York mass \cite{BY} is well defined only for metrics $\g$ for which there is a unique isometric embedding into $\bR^{3}$, 
e.g.~metrics with $K_{\g} > 0$ by the solution to the Weyl problem. Similar convexity requirements are needed so far for 
the Wang-Yau mass \cite{WY}. 

   These and many other definitions of quasi-local mass require solving a complicated system of PDE's to obtain their value; for instance 
solving the isometric immersion problem or the static vacuum Einstein equations. (The main exception is the Hawking mass). 
The behavior of solutions of these equations in the large, non-perturbative realm (e.g.~far from flat) have not yet been well-understood. 

\section{Space-Time Bartnik mass} 

The space-time version of the Bartnik mass, introduced in \cite{Ba1} and \cite{Ba2}, is considerably more complicated than the 
the time-symmetric (static) case \eqref{bm}. The complications are so severe that until very recently, there had 
been essentially no progress on this topic. 

  Let $(\hat M, \hat g, \hat K)$ be a complete AF initial data (Cauchy) surface satisfying the Einstein constraint equations
\be \label{const}
|K|^{2} - H^{2} - R_{g} = \mu, \ \d(K - Hg) = J,
\ee
where $(\mu, J)$ are required to satisfy the dominant energy condition $\mu^{2} \geq |J|^{2}$; we have dropped the hat notation 
for simplicity. The vacuum case $\mu = J = 0$ is of particular interest. 

  As in \S 1, let $(\O, g_{\O}, K_{\O})$ be a smooth bounded domain with boundary (say) $\partial \O = S^{2}$ satisfying 
the constraint equations \eqref{const}, and let $(M, g, K)$ be an AF extension of $\partial \O = \dm$ also satisfying \eqref{const}. 
We view $\O \subset \hat M$, so that $M = \hat M \setminus \O$. Let $\hat \cM$ be the 4-d space-time generated by the 
initial data $(\hat M, \hat g, \hat K)$ by solving the Einstein evolution equations and similarly $\cM \subset \hat \cM$. 

  To begin, one must define suitable boundary data for the quasi-local mass. Here the Bartnik boundary data are given by 
\be \label{stb}
(\g, H, tr_{\dm}K, K(N)^{T}),
\ee
where $N$ is the inward unit normal to $\dm$ in $M \subset \cM$ and $K(N)^{T}$ is the restriction of the 1-form $K(N)$ to $\dm$. 
Note that the boundary conditions $H, K(N)^{T}$ are ``gauge-dependent", i.e. depend on a choice of the space-like hypersurface 
$M$ at $\dm$ or equivalently a choice of time-like unit normal $\nu$ at $\dm$. Such gauge-dependence occurs for many other 
quasi-local masses (e.g.~that of Brown-York \cite{BY} or Wang-Yau \cite{WY}). 

  Bartnik derived the boundary data \eqref{stb} via two distinct but essentially equivalent arguments given in \cite{Ba2}, 
cf.~also \cite{Ba3}. First, the requirement that the energy-momentum terms $(\mu, J)$ are distributionally well-defined 
across the seam or interface $\partial \O = \dm$ of the interior and exterior regions $\O$ and $M$ leads by straightforward 
analysis to the equality of \eqref{stb} at $\partial \Omega$ and $\dm$. The data $(\g, H)$ arise from the 
Hamiltonian constraint, just as in the static case discussed in \S 2, while the data $(tr_{\dm}K, K(N)^{T})$ arise from the 
divergence or momentum constraint in \eqref{const}. 

  A second and in certain respects more significant argument for the choice of the boundary 
conditions \eqref{stb} is that they are exactly the boundary terms arising from the first variation of the Regge-Teitelboim 
Hamiltonian. In more detail, the RT-Hamiltonian is defined by 
\be \label{rt2}
\cH_{RT}(g, K, \xi) = 16\pi \<\xi_{\infty}, \bP_{ADM}\> - \int_{M}\<\xi, \Phi(g,K)\>,
\ee
where $\bP_{ADM}$ is the ADM energy-momentum 4-vector, $\xi = (u, X)$ is the lapse-shift, $\xi_{\infty}$ is their 
value $(u_{\infty}, X_{\infty})$ at space-like infinity and $\Phi$ is the constraint operator $\Phi(g, K) = (\mu, J)$ as in \eqref{const}. 
This is the direct generalization of the static Regge-Teitelboim Hamiltonian \eqref{cL}. 

   One may then derive the first variation of $\cH_{RT}$. Just as in the static case, there are bulk integrals over $M$ and 
boundary integrals at $\dm$. A lengthy but well-known calculation shows that the vanishing of the 
first variation of the bulk integral gives the Euler-Lagrange equations 
\be \label{station}
Ric_{g} - u^{-1}D^{2}u + (tr K)K - 2K^{2} + u^{-1}\cL_{X}K = 0,
\ee
$$\cL_{X}g = 2uK,$$
in addition to the vacuum constraints \eqref{const}. These are a determined system of 16 equations for the 16 unknowns 
$(g, K, u, X)$. As shown by Bartnik in \cite{Ba2}, the vanishing of the variation of the boundary terms gives exactly the 
vanishing of the variation of the boundary data \eqref{stb}. 

  The equations \eqref{station} are the stationary vacuum Einstein equations, giving time-independent solutions of the 
vacuum Einstein equations with Killing field $\partial_{t} = u\nu + X$, where $\nu$ is the time-like unit normal to $M$. 

Writing $\bP_{ADM} = (E, P)$ in the usual way, the ADM mass $m_{ADM}$ is given by 
$$m_{ADM} = \sqrt{|\bP_{ADM}|^{2}} = \sqrt{E^{2} - |P|^{2}}.$$ 

   Each of the definitions of Bartnik mass discussed above in the static case, namely $m_{B}(\O)$ in \eqref{bmor}, $m_{B}(\g,H)$ 
in \eqref{bm}, $m_{B}^{0}$ in \eqref{bm0}, $m_{B}^{out}$ in \eqref{bmout} may be formally extended to the space-time setting. 
For instance, let $\cP^{0}$ be the space 
of triples $(M, g, K)$ with $M \simeq \bR^{3}\setminus B$, with $(g, K)$ asymptotically flat data satisfying the dominant 
energy constraints \eqref{const}, with $(M, g, K)$ containing no apparent horizons surrounding $\dm$. Given Bartnik boundary data 
$(\g, H, tr_{\dm}K, K(N)^{T})$, let 
$$\cP_{(\g,H,tr_{\dm}K,K(N)^{T})}^{0} \subset \cP^{0},$$
be the subset of $\cP^{0}$ with fixed Bartnik boundary data at $\dm$. Then 
\be \label{bmsta}
m_{B}(\g, H, tr_{\dm}K, K(N)^{T}) = \inf \{m_{ADM}(M, g, K): (M, g, K) \in \cP_{(\g,H,tr_{\dm}K,K(N)^{T})}^{0}\}.
\ee

  In  \cite{Ba2}, Bartnik again raised the minimization conjecture that the mass $m_{B}$ in \eqref{bmsta} is realized by a unique AF 
stationary vacuum space-time $(M, g, K)$ with given boundary data \eqref{stb}, and the corresponding stationary extension 
conjecture that all such boundary data are uniquely realized by a stationary vacuum solution. These are of course even 
more difficult than the static versions. All of the issues, problems and conjectures discussed in previous sections apply  
here also, both for the mass \eqref{bmsta} and the analogous space-time version of the modifications of static masses 
discussed previously. 

  However, very little has been known in this space-time setting. As a simple example, the analog of the result of Miao \cite{Mi} 
that static vacuum solutions $(M, g, u)$ have no horizons surrounding $\dm$ is unknown for stationary vacuum solutions. 

  Even the most basic question, raised by Bartnik in \cite{Ba1}, of whether the system \eqref{station} with boundary conditions \eqref{stb} 
is an elliptic system in a suitable gauge was unresolved. This has recently been proved to be true by An \cite{An2}. Given this, 
one has now at least an effective starting point to study some of these conjectures in more detail. 

  Thus, analogous to the static case, let $\bE_{sta}$ be the space of all AF solutions $(g, K)$ of the stationary vacuum Einstein equations 
\eqref{stat} on $M \simeq \bR^{3}\setminus B$ and let $\cE_{sta}$ be the corresponding moduli space, obtained by 
dividing out by a suitable group of diffeomorphisms fixing $\dm$. The choice of the correct diffeomorphism or gauge 
group here is quite subtle, since the group of diffeomorphisms fixing $\dm$ does not preserve the data \eqref{stb} at 
$\dm$. Let then $\cD_{4}$ be the group of time-independent diffeomorphisms of the space-time $\cM$ fixing $\dm$; any element 
in $\cD_{4}$ can be written as a pair $(\psi, f)$, where $\psi \in \cD_{3}$ is a diffeomorphism of $M$ fixing $\dm$ and $f$ 
represents a time translation $f: M \to \bR$ with $f|_{\dm} = 0$, moving $M$ in $\cM$ but fixing $\dm$. The correct choice 
of gauge group $\cD \subset \cD_{4}$, as discussed in \cite{An2}, is 
$$\cD = \{(\psi, f) \in \cD_{4}: N(f) = 0 \ {\rm at} \ \dm\},$$
where $N$ is the unit inward normal to $\dm$ in $M$. This choice is naturally related to the gauge-dependence 
of the boundary conditions \eqref{stb} discussed above. 

  It is then proved in \cite{An1}, \cite{An2} that $\cE_{sta}$ is a smooth infinite dimensional Banach manifold (in a suitable 
$C_{\d}^{m,\a}$ H\"older topology). Moreover, by \cite{An2}, the Bartnik boundary map 
$$\Pi_{B}: \cE_{sta} \to \cB,$$
$$\Pi(g, K) = (\g, H, tr_{\dm}K, K(N)^{T}),$$
is a smooth Fredholm map. In addition, the map $\Pi_{B}$ is a local diffeomorphism in a 
neighborhood of standard, round flat boundary data $(\g_{+1}, 2, 0, 0)$. 

\medskip 

  It would be very interesting to extend the Regge-Teitelboim-Bartnik program discussed in \S 2 to the stationary or space-time 
setting. Thus one would like to study the Hamiltonian \eqref{rt2} and its critical points on the the constraint space 
$\cC_{(\g,H,tr_{\dm}K, K(N)^{T})}$: the space of solutions of the constraint equations \eqref{const} with fixed boundary 
data \eqref{stb}. The first main task in such a program is to prove that $\cC_{(\g,H,tr_{\dm}K, K(N)^{T})}$ can be given 
the structure of a smooth Banach (or Hilbert) manifold. This was proved by Bartnik \cite{Ba3} in the case of empty 
boundary and in the low-regularity regime. Recent progress on this problem has been obtained by An \cite{An3}. 

\bigskip 

{\sl Acknowledgements.} Many thanks to Zhongshan An and Jeff Jauregui for discussions and remarks related to this 
work. 

\bibliographystyle{plain}

\end{document}